\documentclass[12pt,oneside]{amsart}
\usepackage{amsthm,amsfonts,amssymb,euscript}

\renewcommand{\emptyset}{\varnothing}
\newcommand{\mt}[1]{\operatorname{#1}}
\newcommand{\Supp}{\operatorname{Supp}}
\newcommand{\Diff}{\operatorname{Diff}}
\newcommand{\down}[1]{\left\lfloor #1\right\rfloor}
\newcommand{\up}[1]{\left\lceil #1\right\rceil}

\newcommand{\Sing}{\operatorname{Sing}}
\newcommand{\discr}{\operatorname{discr}}
\newcommand{\codim}{\operatorname{codim}}
\newcommand{\Center}{\operatorname{Center}}
\newcommand{\totaldiscr}{\operatorname{totaldiscr}}
\newcommand{\CC}{{\mathbb C}}
\newcommand{\ZZ}{{\mathbb Z}}
\newcommand{\QQ}{{\mathbb Q}}
\newcommand{\PP}{{\mathbb P}}
\newcommand{\NN}{{\mathbb N}}

\newcommand{\OOO}{{\mathcal O}}

\newcommand{\KKK}{{\EuScript{K}}}
\newcommand{\T}{\mathcal{T}}
\newcommand{\Msm}{{\Phi}_{\mt{\mathbf{sm}}}}
\newcommand{\qq}{\mathbin{\sim_{\scriptscriptstyle{\mathbb Q}}}}

\renewcommand{\tilde}[1]{\widetilde{#1}}
\renewcommand{\bar}[1]{\overline{#1}}

\newtheorem{theorem}[subsection]{Theorem}
\newtheorem{conjecture}[subsection]{Conjecture}
\newtheorem{proposition}[subsection]{Proposition}
\newtheorem{lemma}[subsection]{Lemma}
\newtheorem{corollary}[subsection]{Corollary}
\theoremstyle{definition}

\theoremstyle{remark}

\title{On log canonical thresholds, II}
\author{Yu.~G.~Prokhorov}
\thanks{This work was partially supported by the grant
INTAS-OPEN-97-2072}

\address{Department of Algebra, Faculty of Mathematics, 
Moscow State Lomonosov University, 117234 Moscow,
Russia}
\email{prokhoro@mech.math.msu.su}

\address{Current address: Max-Plank-Institut f\"ur Mathematik
Vivatsgasse 7, 53111 Bonn, Germany}
\email{prokhoro@mpim-bonn.mpg.de}

\begin{document}
\begin{abstract}
We prove that the only accumulation points of the set
$\T_3$ of all three-dimensional  
log canonical thresholds in the interval $[1/2,1]$ are 
$1/2+1/n$, where $n\in\ZZ$, $n\ge 3$.
\end{abstract}

\maketitle
\section{Introduction}
In  this paper we continue our study of the structure of
the set $\T_3$ of all three-dimensional  
log canonical thresholds started in \cite{Pr2}.
Notation and results of the Log Minimal 
Model Program \cite{Ut} will be used  freely.
\par
Let $X$ be a normal algebraic variety and let $F$ be an effective
integral non-zero
$\QQ$-Cartier divisor on $X$. Assume that $X$ has at
worst log canonical singularities. \textit{The log canonical
threshold} of $(X,F)$ is defined by
\[
c(X,F)=\sup \left\{c \mid \text{$(X,cF)$
is log canonical}\right\}.
\]

For each $d\in \ZZ$, $d\ge 2$ define the 
following set $\T_d\subset [0,1]$ by
\[
\T_d:= \left\{c(X,F)\ \left|\ 
\text{\begin{tabular}{p{220pt}}
$\dim X=d$, $X$ has only log canonical singularities 
and $F$
is an effective non-zero Weil $\QQ$-Cartier divisor
\end{tabular}}\right.\right\}.
\]
The above does not define $\T_1$ 
but it is naturally to put
\[
\T_1:=\left\{\left. \frac1n\quad \right| \quad 
n\in \NN\cup\{\infty\}\right\}.
\]

The sets $\T_d$ have rather  inductive nature:
it is easy to show that $\T_{d-1}\subset \T_d$
and $\partial\T_d\supset \T_{d-1}$ (see \cite[8.21]{Ko}), 
where $\partial\T$  
is the set of all accumulation points
of $\T$. 

\begin{conjecture}[{\cite{Ko}}]
\label{conj}
The accumulation set $\partial \T_d$ of 
$\T_d$ is precisely $\T_{d-1}$. 
\end{conjecture}
This conjecture is the only one instance where 
the such a phenomena occurs.
The similar behavior is expected for the fractional indices
of log Fano varieties \cite{Sh0}, \cite{A-1}, 
minimal log discrepancies \cite{Sh0}, \cite{Sh1}, \cite{B},
Kodaira energy \cite{F} etc.

In dimension two Conjecture \ref{conj} easily follows 
from explicit description of $\T_2$ \cite{Ku1}.
In this paper we generalize the result of \cite{Pr2}
and prove Conjecture \ref{conj}  in dimension three
for the interval $\left[\frac12,1\right]$:
 
\begin{theorem}
\label{main}
\begin{equation}
\label{eq-Al}
\partial \T_3\cap \left[\frac12,1\right]=\T_{2}\cap 
\left[\frac12,1\right]
=\left\{\left.\frac12+\frac1n\quad \right|\quad
n\in\ZZ, n\ge 3\right\}. 
\end{equation}
\end{theorem}
Note that \eqref{eq-Al} very similar to 
the corresponding results for log Del Pezzo surfaces 
\cite{A-1}. 
Our proof based on inductive arguments and 
boundedness result \cite{A}. As an intermediate result, we 
prove the following
easy but very important fact:

\begin{proposition}
\label{main1}
Assume the LMMP in dimension $d$.
Let $X\ni o$ be a $d$-dimensional $\QQ$-factorial
log terminal singularity\footnote{By \cite[Lemma 4.1]{Pr2} 
computing $\T_d$ we can consider
only those singularities $X$ which are $\QQ$-factorial and 
log terminal.} and let 
$F$ be an (integral) Weil divisor on $X$.
Let $c:=c(X,F)$ be the log canonical threshold.
Then one of the following holds:
\begin{enumerate}
\item
$c\in \T_{d-1}$; or
\item
$c\notin \T_{d-1}$ and there is exactly one divisor 
$S$ of the function field 
$\KKK(X)$ with discrepancy $a(S,cF)=-1$ (i.e., the pair 
$(X,cF)$ is \emph{exceptional} in the sense of \cite{Sh}).
\end{enumerate}
Moreover in case {\rm (ii)}, $\Center(S)=o$. 
\end{proposition}

\subsubsection*{Acknowledgments}
This work was carried out during my stay at 
Max-Planck-Institut f\"ur Mathematik.
I would like to thank MPIM for wonderful working environment.

\section{Preliminary results}
\subsection*{Notation}
All varieties are assumed to be algebraic varieties defined over
the field $\CC$. A \textit{log variety} (or a \textit{log pair}) $(X,D)$
is a normal quasiprojective variety $X$ equipped with a \textit{boundary}
that is
a $\QQ$-divisor $D=\sum d_iD_i$ such that $0\le d_i\le 1$ for all $i$.
We use terminology, definitions and abbreviations of the Log Minimal Model 
Program \cite{Ut}. Recall that $a(E,D)$ denotes the discrepancy of 
$E$ with respect to $D$ and 
\[
\begin{array}{lll}
\discr(X,D)&=&\inf_E\{ a(E,D)\mid \codim\Center(E)\ge 2\}.
\vspace{6pt}\\
\totaldiscr(X,D)&=&\inf_E\{ a(E,D)\mid \codim\Center(E)\ge 1\}.
\end{array}
\]
Recall also our notation of \cite{Pr2}:
\[
\begin{array}{lll}
\Msm&=&\bigl\{1-\frac1m \mid m\in\NN\cup\{\infty\}\bigr\},
\vspace{6pt}\\
\Msm^\alpha&=&\Msm\cup [\alpha,1],\quad\text{for $\alpha\in[0,1]$}.
\end{array}
\]
Let $\Phi$ be any subset of $\QQ$ and let $D=\sum D_i$
be a $\QQ$-divisor. We write $D\in\Phi$ if $d_i\in \Phi$ for all $i$.

\begin{lemma}
\label{P1}
Fix a constant $N\in\ZZ$, $N\ge 6$. 
Let $\Lambda=\sum_{i=1}^r\lambda_i\Lambda_i$ 
be a boundary on $\PP^1$ such that 
\begin{enumerate}
\item
$K_{\PP^1}+\Lambda\equiv 0$;
\item
$\Lambda\in\Msm^{\frac12+\frac1N}$; and
\item
$1> \lambda_j>\frac12+\frac1{N}$ for some $j$.
\end{enumerate}
Then $\lambda_i\le 1-\frac1{N}$ for all $i$.
\end{lemma}

\begin{proof}
Clearly, $r=3$ and $\down{\Lambda}=0$. 
Assume that $\lambda_1>1-\frac1N$.
Then 
\[
1<\lambda_2+\lambda_3<1+\frac1{N}.
\]
Since $\lambda_i\ge \frac12$, we have
$\lambda_2,\lambda_3<\frac12+\frac1{N}$. Thus 
$\lambda_2=\lambda_3=\frac12$ and $\lambda_1=1$,
a contradiction.
\end{proof}

\begin{lemma}
\label{LCT-2}
Fix a constant $N\in\NN$, $N\ge 6$.
Let $(S\ni o,\Theta=\sum\vartheta_i\Theta_i)$ 
be a klt log surface germ with $\Theta\in\Msm^{\frac12+\frac1N}$.
Define the following boundary 
$\Xi$ with $\Supp(\Xi)=\Supp(\Theta)$:
\begin{equation}
\label{eq-Xi1}
\Xi:=\sum \xi_i\Theta_i,\qquad \xi_i=
\begin{cases}
1 & \text{if $\vartheta_i> 1-\frac1{N}$}, \\ \vartheta_i &
\text{otherwise}.
\end{cases}
\end{equation}
Then $(S,\Xi)$ is lc.
\end{lemma}

\begin{proof}
If $\vartheta_i\le 1-\frac1{N}$ for all $i$, there is nothing to prove.
Assume that $\Xi\neq \Theta$ and
$(S,\Xi)$ is not lc. Replacing $\Theta$ with
$\Theta+\alpha(\Xi-\Theta)$, $\alpha>0$, 
we may assume that 
$(S,\Theta)$ is lc but not klt (and $\down{\Theta}=0$).
Let $\mu\colon \bar S\to S$ be an 
\emph{inductive blowup}\footnote{In \cite{Pr1}
such a $\mu$ was called 
\emph{plt-blowup of the pair $(S,\Theta)$}.}
of the pair $(S,\Theta)$ (see \cite[Prop. 5]{Pr1})
and let $E$ be the exceptional divisor.  
By definition, $E$ is irreducible,  $a(E,\Theta)=-1$ and
$(\bar S,E)$ is plt. 
Write 
\[
\mu^*(K_S+\Theta)=K_{\bar S}+E+\bar\Theta,
\]
where $\bar\Theta$ is the proper transform of $\Theta$.
Clearly, $\mu(E)\in \Theta_j$ with $\vartheta_j>1-\frac1N$.

\subsection{}
\label{new}
By \cite[Corollary~2.5]{Pr2}, $\Diff_E(\bar \Theta)
\in\Msm^{\frac12+\frac1N}$. 
Pick a point $\bar P\in E\cap \bar \Theta_j$. 
Then $\bar\Theta_j$ is the only component of $\bar\Theta$,
passing through $\bar P$ (see \cite[Corollary~2.4]{Pr2}).
Moreover, 
$(\bar S,E+\bar\Theta_j)$ is lc at $\bar P$
\cite[Lemma 3.2]{Pr2}. Hence
$(\bar S,E+\bar\Theta)$ is
plt at $\bar P$ and the coefficient
$\lambda'$ of $\Diff_E(\bar \Theta)$
at $\bar P$ satisfies the inequality
$1-\frac1N<\lambda'<1$. Therefore,
$\Lambda:=\Diff_E(\bar \Theta)$ satisfies conditions of 
Lemma~\ref{P1}. This gives us
$\Diff_E(\bar \Theta)\in \left[1,\frac1N\right]$, 
a contradiction.
\end{proof}

\begin{lemma}
\label{pair-discr}
Let $(S\ni o,\Lambda=\sum \lambda_i\Lambda_i)$ be a
log surface germ such that 
$\Lambda\in \left(1-\frac1{N},1\right]$.
Assume
that $\discr(S,\Lambda)\ge -1+\frac1{N}$ at $o$ for 
$N\in\ZZ$, $N\ge 6$. Then
$\sum \lambda_i\le 2-\frac1{N}$. In particular, $\Lambda$ has at 
most two components.
\end{lemma}

\begin{proof}
For some $\Lambda':=\Lambda+t(\up{\Lambda}-\Lambda)$,
$0<t\le 1$
the pair $(S,\Lambda')$ 
is lc but not plt at $o$. By Lemma~\ref{LCT-2}, we have
$\Lambda'=\up{\Lambda}$, i.e., $(S,\up{\Lambda})$ 
is lc. 
If $\Lambda$ has only one component, there is nothing 
to prove. So, we may assume that $\Lambda$ has exactly 
two components \cite[Ch. 3]{Ut}.
Then near $o$ we have
\[
(S,\up{\Lambda})\simeq_{\mt{an}} (\CC^2,\{xy=0\},0)/\ZZ_m(1,q),
\]
where $m\in\NN$ and $\gcd(m,q)=1$.
Take $q$ so that $1\le q\le m$. As in
the proof of Lemma 3.3 in \cite{Pr2}, 
considering the weighted blow up with weights $\frac1m(1,q)$ we
get $\lambda_1+\lambda_2\le 2-\frac1{N}$.
\end{proof}

\section{Proof of Proposition~\ref{main1}. Corollaries}
Notation and assumption as in  Proposition~\ref{main1}.
Let $f\colon Y\to X$ be an inductive blowup of the pair $(X,cF)$
(see \cite[Prop. 5]{Pr1}). 
Write 
\[
f^*(K_X+cF)=K_Y+cF_Y+S,
\] 
where $F_Y$ is the proper
transform of $F$ and $S$ is the (irreducible) exceptional divisor. 
By definition, $(Y,S)$ is plt.

Assume that $c\notin\T_{d-1}$.
If $f(S)\neq o$, then the pair $(X,cF)$ is lc but not
klt along $f(S)$. 
Taking the general hyperplane section we get $c\in \T_{d-1}$. 

Hence $f(S)=o$. It is sufficient to show that $(Y,S+cF_Y)$ is plt
(see \cite[3.10]{Ko}). 
Assume the converse. Then there is an divisor $E\neq S$ 
of the field $\KKK(Y)$ such
that $a(E,S+cF_Y)=-1$. Since $(Y,S)$ is plt,
$\Center_Y(E)\subset E\cap F_Y$.

Pick a point $P\in \Center_Y(E)$ and consider $Y$ 
as a germ near $P$.
Take the minimal
$m\in \NN$ such that $mS\sim 0$ near $P$ and
let 
\[
Y':=\mt{Spec}\Bigl(\bigoplus_{i=0}^{m-1}
\OOO_Y(iS)\Bigr).
\]
Then the projection $\varphi\colon Y'\to Y$ is
an \'etale in codimension one 
$\ZZ_m$-covering. 
Put $P':=\varphi^{-1}(P)$,
$F_Y':=\varphi^{*}F_Y$, and $S':=\varphi^{*}S$.
Then $(Y',S')$ is plt and $(Y,S'+cF_Y')$ is lc but not plt
near $P'$ (see \cite[\S 2]{Sh}).
Since $S'$ is Cartier, $\Diff_{S'}(0)=0$
(i.e., no codimension two components of 
$\Sing(Y')$ are contained in $S'$).
By the Adjunction \cite[Th. 17.6, 17.7]{Ut}
$(S',cF_Y'|_{S'})$ is lc but not klt near $P'$.
Hence $c=c(S',F_Y'|_{S'})$ and $c\in\T_{d-1}$.

The Adjunction \cite[17.6]{Ut} and \cite[Cor. 3.10]{Sh} 
gives us the following:
\begin{corollary}
\label{description}
Let $c\in \T_d\setminus \T_{d-1}$. 
Assume that the LMMP in dimension $d$ holds.
Then there is a log pair $(S,\Theta)$ such that 
\begin{enumerate}
\item
$(S,\Theta)$ is klt;
\item
$K_S+\Theta\qq 0$;
\item
$\Theta=\sum_i\vartheta_i\Theta_i$, where 
\[
\vartheta_i=1-\frac1{m_i}+\frac{k_ic}{m_i},\quad 
m_i\in\NN,\ k_i\in\ZZ_{\ge 0},\ k_ic<1;
\]
\item
$-\left(K_S+\sum_i(1-\frac1{m_i})\Theta_i\right)$ is ample. 
In particular, $\sum k_i>0$.
\end{enumerate}
\end{corollary}

\begin{corollary}
\label{T2}
Let $c\in \T_2\setminus \T_1$. 
Then there are $m_i\in\NN$, $k_i\in\ZZ_{\ge 0}$
such that 
\begin{equation}
\label{eq-m}
\text{$k_ic<1$, $\sum k_i>0$, and}\quad
\sum_i\left(1-\frac1{m_i}+\frac{k_ic}{m_i}\right)=2.
\end{equation}
Moreover,
allowing $k_ic=1$ in \eqref{eq-m}, we get 
$c=\frac1{k_i}\in \T_1\subset\T_2$.
Conversely, 
if there are 
$m_i\in\NN$, $k_i\in\ZZ_{\ge 0}$
satisfying \eqref{eq-m}, then $c\in\T_2$.
\end{corollary}

\begin{proof}
Apply Corollary~\ref{description}. We obtain
$S\simeq\PP^1$ and $\deg \Theta=2$. 
The inverse implication follows by \cite{Ku1}. 
\end{proof}

\begin{corollary}[\cite{Ku1}]
Any $c\in\T_2\cap (\frac12,1]$ has the following form
\[
\frac12+\frac1n,\quad n\in\ZZ,\ n\ge 2.
\]
\end{corollary}

\subsection{}
For $c\in [0,1]\cap \QQ$,
let $\mathcal{LP}(c)$ be the class of all projective
klt log surfaces $(S,\Theta)$ 
satisfying conditions {\rm (i)-(iv)}
of Corollary~\ref{description}. 
Then 
\[
\T_3\setminus \T_2\subset \{c \mid \mathcal{LP}(c)\neq\emptyset\}.
\]

\begin{lemma}
\label{curve}
Let $c$ and $(S,\Theta)$ be as in Corollary~\ref{description}
with $d=3$.
Assume that there is a contraction 
$g\colon S\to W$ onto a curve. Then 
all components $\Theta_i$ with $k_i>0$ are 
vertical (i.e., $g(\Theta_i)\neq W$).
\end{lemma}
\begin{proof}
Assume that there is a horizontal component 
$\Theta_i$ with $k_i>0$. Let $S_w$ be the general fiber.
Then $S_w\simeq\PP^1$ and by Adjunction we have 
equality \eqref{eq-m}:
\[
\deg \Theta|_{S_w}=
\sum_{\Theta_i\cap S_w\neq\emptyset}
\left(1-\frac1{m_i}+\frac{k_ic}{m_i}\right)=2.
\]
By our assumption, $\sum_{\Theta_i\cap S_w\neq\emptyset} k_i>0$.
Thus $c\in \T_2$, a contradiction.
\end{proof}

\begin{corollary}
\label{description1}
Let $c\in \T_3\setminus \T_{2}$. 
Then there is a log surface $(S,\Theta)\in \mathcal{LP}(c)$ with
$\rho(S)=1$.
\end{corollary}
\begin{proof}
Denote 
\[
\Theta^c:=\sum_{k_i>0} \left(1-\frac1{m_i}+
\frac{k_ic}{m_i}\right)\Theta_i
\]
and run 
$K_S+\Theta-\Theta^c$-MMP.
Since $K_S+\Theta\equiv 0$,
each time we contract an extremal ray $R$ such that 
$R\cdot \Theta^c>0$. Hence $\Theta^c$ is not contracted.  
By Lemma~\ref{curve}, at the end we obtain a model with $\rho=1$. 
\end{proof}

\section{Proof of the main theorem}
\label{proof}
In this section we prove Theorem \ref {main}.

\begin{lemma}
\label{bound}
For any $\epsilon>0$ and $\frac12>\xi>0$ there exists a finite set 
$\mathcal{M}_{\xi,\epsilon}\subset [0,1]$ such that $c\in
\mathcal{M}_{\xi,\epsilon}$ whenever $c>\xi$ and there is 
$(S,\Theta)\in \mathcal{LP}(c)$ with 
\[
\totaldiscr(S,\Theta)>-1+\epsilon.
\]
\end{lemma}
\begin{proof}
Since $\Theta\neq 0$, one can apply \cite[Th. 6.9]{A} 
to $(S,\Theta)$. This gives us that the family
$\mathbf{S}$ of all such $S$ is bounded.
That is there is a family $\mathbf{S}\to \mathbf{H}$
such that every $S$ is a fiber of $\mathbf{S}\to \mathbf{H}$.
Therefore there is a polarization $\mathbf{L}$ on $\mathbf{S}$ 
giving us an embedding $\mathbf{S}\hookrightarrow \PP$ 
over $\mathbf{H}$. This induces 
a very ample divisor $L$ on each $S\in \mathbf{S}$.
For all coefficients of $\Theta$ we have 
$\vartheta_i> \xi$. 
Then 
\[
L\cdot \sum_i \Theta_i< -\frac1{\xi} L\cdot K_S\le 
\mt{Const}(\epsilon,\xi). 
\]
Hence the family
of all $\sum \Theta_i$ is represented by a closed subscheme
of $\PP$ over $\mathbf{H}$.
This shows that the pair $(S,\Supp(\Theta))$ is bounded.
From the equality
\[
L\cdot K_S+L\cdot \Theta=0
\]
we obtain the following linear equation in $c$:
\[
L\cdot K_S+\sum_i \left(1-\frac1{m_i}+\frac{k_ic}{m_i}\right) 
(L\cdot\Theta_i)=0,
\]
where 
\[
1-\frac1{m_i}+\frac{k_ic}{m_i}\le -\totaldiscr(S,\Theta)<1-\epsilon.
\]
This gives us a finite number of possibilities for the $m_i$, $k_i$ and
$c$.
\end{proof}

\begin{lemma}
\label{2}
Fix constants $N\in\ZZ$, $N\ge 6$ and $0<\epsilon<\frac1{N}$.
Let $(S,\Theta=\sum \vartheta_i\Theta_i)$
be a klt log surface such that $\Theta\in\Msm^{\frac12+\frac1N}$.
Assume that there are at least two divisors 
of the function field $\KKK(S)$ such that 
\[
a(\phantom{E_i},\Theta)<-1+\frac1{N}-\epsilon.
\]
Then 
\[
\totaldiscr(S,\Theta)> -1+\epsilon.
\]
\end{lemma}
\begin{proof}
Let 
$\mu\colon \tilde S\to S$ be the blowup 
of all divisors with discrepancies
$a(\phantom{E_1},\Theta)<-1+\frac1{N}$
(see \cite[Th. 17.10, 2.12.2]{Ut}) and let $\tilde \Theta$
be the crepant pullback of $\Theta$:
\[
K_{\tilde S}+\tilde \Theta=\mu^*(K_S+\Theta),
\quad \mu_*\tilde \Theta=\Theta.
\] 
Then
$(\tilde S,\tilde \Theta)$ satisfies conditions of Lemma~\ref{2}.
Moreover,
\[
\discr(\tilde S,\tilde\Theta)\ge -1+\frac1{N}.
\]
Clearly,
\[
\totaldiscr(S,\Theta)=
\totaldiscr(\tilde S,\tilde \Theta)
\]
(see \cite[3.10]{Ko}).
Replace $(S,\Theta)$ with $(\tilde S,\tilde \Theta)$.
Up to permutations of the $\Theta_i$ we may assume that
\[
\vartheta_1,\vartheta_2>1-\frac1{N}+\epsilon. 
\]
Now it is sufficient to show that $\vartheta_i<1-\epsilon$
for all $i$. 
Consider the boundary 
$\Xi$ with $\Supp(\Xi)=\Supp(\Theta)$ as in \eqref{eq-Xi1}.
Then $\down{\Xi}=\up{\Xi-\Theta}$. 
For a sufficiently small positive rational $\alpha$, the $\QQ$-divisor
$\Theta -\alpha (\Xi -\Theta)$
is a boundary. 
It is clear
that 
\[
K_{S}+\Theta -\alpha (\Xi -
\Theta)\equiv -\alpha (\Xi -
\Theta)
\]
cannot be nef. 
By Lemma~\ref{LCT-2} the pair $(S, \Xi)$ is lc.

Run $K_{S}+ \Theta -\alpha (\Xi -\Theta)$-MMP. 
On each step we contract an extremal ray 
$R$ such that 
\[
(K_{S}+\Xi)\cdot R=(\Xi -\Theta)\cdot R>0.
\] 

\subsection{}
We claim that none of the components of $\down{\Xi}$ is contracted.
Indeed, assume that $\varphi\colon S\to S^o$ contracts 
$C\subset \down{\Xi}$. Take $\Theta':=\Theta+\beta C$
so that $\down{\Theta'}=C$ and $\Theta'\le \Xi$.
Since  $(K_S+\Xi)\cdot C>0$ and $(K_S+\Theta)\cdot C<0$, 
there is a component,
say $\Theta_0$, of $\down{\Xi}$ meeting $C$.
Further, take $\Theta'':=\Theta'+\gamma(\Xi-\Theta')$ so that 
$(K_S+\Theta'')\cdot C=0$. Then $0<\gamma<1$. 
It is easy to see that $\Theta''\in\Msm^{\frac12+\frac1N}$
and $\down{\Theta''}=C$. Note that $(S,\Theta'')$ is lc
(because so is $(S,\Xi)$).
As in the \ref{new}, we can 
apply  Lemma~\ref{P1} to $\Diff_C(\Theta''-C)$ to derive 
a contradiction. This  proves our claim.

\subsection{}
By Lemma~\ref{LCT-2} the lc property of $(S, \Xi)$ 
is preserved on each step.
At the end of the MMP we get
a birational model $(\bar S, \bar\Theta)$ 
with nonbirational extremal 
$\bar\Xi-\bar\Theta$-positive  contraction 
$g\colon \bar S\to W$, where $W$ is either a curve
or a point.

\subsubsection{Subcase: $W$ is a curve}
\label{c}
Then $\rho(\bar S)=2$.
Let $\bar S_w$ be the general fiber of $g$.  
Then $\Diff_{\bar S_w}(\Theta)$ satisfy conditions of Lemma~\ref{P1}.
This yields a contradiction.

\subsubsection{Subcase: $W$ is a point}
Then $\rho(\bar S)=1$ and every two components of 
$\Theta$ intersects each other. By Lemma~\ref{pair-discr},
\[
\vartheta_1\le 2-\frac1{N}-\vartheta_2<1-\epsilon.
\]
Similarly, if  $i\neq 1$ and the image of $\Theta_i$
on $\bar S$ is not a point, then
\[
\vartheta_i\le 2-\frac1{N}-\vartheta_1<1-\epsilon.
\]
But if $\Theta_i$ is contracted to a point  
on $\bar S$, then $\Theta\not\subset\down{\Xi}$.
In this case, $\theta_i\le 1-\frac1N<1-\epsilon$.
This proves our lemma.
\end{proof}

\subsection{}
Now we are ready to prove Theorem~\ref{main}.
Assume that there is a sequence $c_n\in \T_3\cap [\frac12,1]$
such that $c_{n_1}\neq c_{n_2}$ for $n_1\neq n_2$
and $\lim c_n=c_{\infty}\notin\T_2$.
Take constants $N\in\NN$ and $\epsilon\in\QQ$ so that 
\[
\begin{array}{l}
N\ge 6,\qquad \frac12+\frac1{N}<c_{\infty}, \ \text{and} 
\vspace{8pt}\\
0<\epsilon<\min\left\{c_{\infty}-\frac12-\frac1N,\ \frac1N\right\}.
\end{array}
\]

By passing to a subsequence, we may assume that
$c_n>\frac12+\frac1{N}+\epsilon$ for all $n$. 
For every $c_n$ we have the
corresponding log surface 
$(S_n,\Theta_n)\in\mathcal{LP}(c_n)$ 
with $\rho(S_n)=1$
(see Corollaries~\ref{description} and \ref{description1}).
In particular, $\Theta_n\in \Msm^{\frac12+\frac1N+\epsilon}$.
Write $\Theta_n=\sum_i\vartheta_{n,i}\Theta_{n,i}$.
By construction, 
\begin{equation}
\label{eq-va1}
\vartheta_{n,i}=
1-\frac1{m_{n,i}}+\frac{k_{n,i}c_n}{m_{n,i}},
\quad k_{n,i}c_n<1, \quad \sum_i k_{n,i}>0.
\end{equation}

If 
\[
\lim_{n\to \infty} \totaldiscr(S_n,\Theta_n)>-1,  
\]
we can take $\nu>0$ so that
$\totaldiscr(S_n,\Theta_n)\ge-1+\nu$
for $n\gg 0$, 
then $c_n$ belongs to a finite  set 
$\mathcal{M}_{\frac12+\frac1N,\nu}$ by Lemma~\ref{bound}.
This contradicts to our choice of the sequence $c_n$.
From now on we assume that 
\begin{equation}
\label{eq-lim}
\lim_{n\to \infty} \totaldiscr(S_n,\Theta_n)=-1,  
\end{equation}
In particular, 
\[
\totaldiscr(S_n,\Theta_n)< -1+\frac1{N}-\epsilon\
\text{for all $n$.}
\]

Assume that for $n\gg 0$ there are at least two divisors 
of the field $\KKK(S_n)$ with discrepancies
$a(\phantom{E_i},\Theta_n)<-1+\frac1{N}-\epsilon$. 
Then
$(S_n,\Theta_n)$ satisfies conditions of Lemma~\ref{2}.
Therefore 
\[
\totaldiscr(S_n,\Theta_n)>-1+\epsilon,
\]
This contradicts \eqref{eq-lim}.

\subsection{Main case}
Finally we consider the case when for $n\gg 0$ 
there is exactly one divisor 
$\Gamma_n$ with 
\[
\gamma_n:=-a(\Gamma_n,\Theta_n)>1-\frac1{N}+\epsilon.
\] 
We construct a new birational model
$(\bar S_n,\gamma_n\bar\Gamma_n+\bar\Theta_n)$
of $(S_n,\Theta_n)$ with $\rho(\bar S_n)=1$
and such that the center of  $\Gamma_n$ on $\bar S_n$
is a curve. 

\subsubsection{}
\label{c1}
If $\Center_{S_n}(\Gamma_n)$ is a curve, 
then $\Gamma_n=\Theta_{n,i}$ and 
$\gamma_n=\vartheta_{n,i}$ for some $i$.
In this case we just put 
$\bar S_n:=S_n$ and
$\bar \Theta_n:=\Theta_n-\gamma_n\Gamma_n$.
Thus
\[
\bar\Theta_n=\sum_i \bar\vartheta_{n,i}
\bar \Theta_{n,i},
\]
where $\bar \Theta_{n,i}:=\Theta_{n,i}$ whenever 
$\Theta_{n,i}\neq \Gamma_n$ and
\[
\bar \vartheta_{n,i} =
\begin{cases}
0&\text{if $\Theta_{n,i}= \Gamma_n$,}\\
\vartheta_{n,i}
&\text{otherwise.}
\end{cases}
\]

\subsubsection{}
\label{c2}
If  $\Center_{S_n}(\Gamma_n)$ is a point, we consider
the blowup of this $\Gamma_n$:
$\mu\colon \tilde S_n\to S_n$ \cite[Th. 17.10]{Ut}.
Clearly, $\rho(\tilde S_n)=2$.
Write
\[
\begin{array}{l}
K_{\tilde S_n}+\gamma_n\Gamma_n+
\tilde \Theta_n=\mu^*\left(K_{S_n}+\Theta_n\right),
\vspace{8pt}\\
\tilde \Theta_n=\sum \vartheta_i\tilde \Theta_{n,i},\quad
\text{where $\mu_*\tilde \Theta_{n,i}=\Theta_{n,i}$.}
\end{array}
\]

By construction, 
$\vartheta_{n,i}\le 1-\frac1{N}+\epsilon$.
The divisor
$K_{\tilde S_n}+\tilde\Theta_n\equiv-\gamma_n\Gamma_n$
cannot be nef. Therefore, there is a $\Gamma_n$-positive
extremal contraction $\varphi\colon \tilde S_n\to \bar S_n$, 
where $\rho(\bar S_n)=1$. 
By Lemma~\ref{LCT-2}, $(\tilde S_n,\Gamma_n+
\tilde \Theta_n)$ is lc. If $\bar S_n$ is a curve,
we derive a contradiction as in \ref{c}.

Therefore $\varphi$ is birational.
Put $\bar\Theta_n:=\varphi_*\tilde \Theta_n$,
$\bar\Theta_{n,i}:=\varphi_*\tilde \Theta_{n,i}$,
and $\bar\Gamma_n:=\varphi_*\Gamma_n$.
Then $(\bar S_n,\gamma_n\bar\Gamma_n+\bar\Theta_n)$
is klt and $K_{\bar S_n}+\gamma_n\bar\Gamma_n+\bar\Theta_n$
is numerically trivial. Again by Lemma~\ref{LCT-2}, 
$(\bar S_n,\bar\Gamma_n+\bar\Theta_n)$
is lc.

Further, 
\[
\bar\Theta_n=\sum_i \bar\vartheta_{n,i}
\bar \Theta_{n,i},
\]
where 
\[
\bar \vartheta_{n,i} =
\begin{cases}
0&\text{if $\varphi(\tilde \Theta_{n,i})$ is a point,}\\
\vartheta_{n,i}
&\text{otherwise.}
\end{cases}
\]

\subsubsection{}
In both cases \ref{c1} and \ref{c2} we have 
\begin{equation}
\label{eq-va}
\bar \vartheta_{n,i}\le 1-\frac1{N}+\epsilon.
\end{equation}

As in the proof of Lemma~\ref{bound}, apply 
\cite[Th. 6.9]{A} to $(\bar S_n,\bar\Theta_n)$.
We get that the family of all 
$(\bar S_n,\Supp(\bar\Theta_n+\bar\Gamma_n))$
is bounded. By passing to a subsequence we may assume that 
all the discrete invariants $(\bar \Gamma_n)^2$, 
$\bar \Gamma_n\cdot K_{\bar S_n}$, 
$\Theta_{n,i}\cdot K_{\bar S_n}$, $p_a(\bar \Gamma_n)$, and
$K_{\bar S_n}^2$  do no depend on $n$. For short denote them 
by $\bar \Gamma^2$, $\bar \Gamma\cdot K_{\bar S}$, 
$\Theta_i\cdot K_{\bar S}$, $p_a(\bar \Gamma)$, and
$K_{\bar S}^2$, respectively.

From \eqref{eq-va} by
passing to a subsequence we may assume that all constants 
$m_{n,i}$ and $k_{n,i}$ in \eqref{eq-va1} also do not depend on $n$:
\[
\bar \vartheta_{n,i} =1-\frac1{m_i}+
\frac{k_ic_n}{m_i}.
\]
By the Adjunction \cite[Ch. 16]{Ut},
\[
K_{\bar\Gamma_n}+\Diff_{\bar\Gamma_n}\left(
\bar \Theta_n \right)\equiv
(1-\gamma_n)\bar \Gamma_n|_{\bar\Gamma_n},
\]
where $\Diff_{\bar\Gamma_n}\left(
\bar \Theta_n \right)\ge 0$. 
Since $(\bar S_n,\bar\Gamma_n+\bar\Theta_n)$
is lc, $\Diff_{\bar\Gamma_n}\left(
\bar \Theta_n \right)$ is a boundary
(see \cite[Prop. 16.6]{Ut}). 
The coefficients of $\Diff_{\bar\Gamma_n}\left(
\Theta_n \right)$ have the same form as the coefficients 
of $\Theta_n$:
\[
\Diff_{\bar\Gamma_n}\left(
\Theta_n \right)=
\sum_j \left(1-\frac1{s_j}+
\frac{r_jc_n}{s_j}\right)P_j,
\]
where $n_j\in\NN$, $r_j\in\ZZ_{\ge 0}$,
and $r_jc_n\le 1$ (see \cite[Lemma 4.2]{Sh}).
Thus 
\begin{equation}
\label{eq-s}
\sum_j \left(1-\frac1{s_j}+
\frac{r_jc_n}{s_j}\right)
=2- 2p_a(\bar \Gamma)+
(1-\gamma_n) \bar\Gamma^2.
\end{equation}
Here $\bar\Gamma^2>0$, $1-\frac1N+\epsilon<\gamma_n<1$ and 
$p_a(\bar\Gamma)\in\ZZ_{\ge 0}$.
If $r_j=0$ for all $j$, then $\gamma_n$ can be found from
the equation 
\[
\sum_j \left(1-\frac1{s_j}\right)=2- 2p_a(\bar \Gamma)+
(1-\gamma_n) \bar\Gamma^2.
\]
In this case, $\gamma:=\gamma_n$ does not depend on $n$ and
$\gamma<1$. Therefore, 
\[
\totaldiscr(S^n,\Theta^n)>-\gamma>-1.
\]
This contradicts our assumption \eqref{eq-lim}.

Assume that there is at least one component with $r_i=1$.
Passing to the limit as $n\to \infty$ in \eqref{eq-s} we obtain
\[
\sum_j \left(1-\frac1{s_j}+
\frac{r_jc_{\infty}}{s_j}\right)=2- 2p_a(\bar \Gamma)+
(1-\gamma_{\infty}) \bar\Gamma^2.
\]
If $\gamma_{\infty}<1$, then 
\[
\lim_{n\to\infty}\totaldiscr(S^n,\Theta^n)
\ge \min\left\{-\gamma_{\infty},\quad -1+\frac1N-\epsilon\right\}
>-1.
\]
Again we have a contradiction with \eqref{eq-lim}.
Hence $\gamma_{\infty}=1$ and 
\[
0<\sum_j \left(1-\frac1{s_j}+
\frac{r_jc_{\infty}}{s_j}\right)=2- 2p_a(\bar \Gamma).
\]
This gives us that $p_a(\bar \Gamma)$.
By Lemma~\ref{description},
$c_{\infty}\in\T_2$.
Theorem~\ref{main} is proved.

\end{document}